\documentclass{article}
\usepackage{hyperref,multicol,amsmath,amsthm,graphtex3,enumerate}
\usepackage{geometry}
\title{Score sequences of bitournaments}
\author{Severino V.~Gervacio\\
Mathematics and Statistics Department\\
De La Salle University-Manila\\
2401 Taft Avenue, 0922 Manila, Philippines}
\date{}
\newtheorem{thm}{Theorem}[section]
\newtheorem{thmcor}{Corollary}[thm]

\theoremstyle{definition}
\newtheorem{defn}{Definition}[section]
\newcommand\sequence[2]{\ensuremath\langle #1_i\rangle_{i=1}^#2}
\newtheorem{exmp}{Example}[section]
\begin{document}
	\maketitle
	\thispagestyle{empty}
	\begin{abstract}
		The score of a vertex $x$ in an oriented graph is defined to be its outdegree, \emph{i.e.}, the number of arcs with initial vertex $x$. The score sequence of an oriented graph is the sequence of all scores arranged in nondecreasing order.
		
		An oriented complete bipartite graph is called a bitournament.  The score sequence of a bitournament consists of two nondecreasing sequences of nonnegative integers, one for each of the two partite sets.
		
		Moon has characterized the score sequences of bitournaments.  This paper introduces the concept of trimming a sequence and gives a characterization of score sequences of bitournaments utilizing this concept.
	\end{abstract}
\section{Introduction}
We consider here only finite graphs having no loops nor multiple edges. An edge $\{x,y\}$ of a graph $G$  can be changed to the directed edge $xy$ or $yx$.  If every edge of the graph $G$ is changed to a directed edge, the result is called an \emph{oriented} graph.

The \emph{outdegree} of a vertex $x$ in an oriented graph is the total number of directed edges of the form $xy$, and is denoted by $d^+(x)$.  The \emph{indegree} of the vertex $x$, denoted by $d^-(x)$ is the total number of directed edges of the form $yx$.

An oriented complete graph is called a \emph{tournament}.  The score of a vertex $x$ in a tournament is defined as the outdegree of the vertex.  The sequence of scores in non-decreasing order is called the \emph{score sequence} \cite{gervacio} of the tournament.

Landau \cite{landau} characterized the score sequences of tournaments in 1953.

\begin{thm}[\rm Landau]
	The nondecreasing sequence $\sequence s n$ of nonnegative integers $s_i$ is the score sequence of a tournament of order $n$ if and only if
	$$\sum_{i=1}^k s_i\ge \binom k 2$$
	for $k=1, 2, \ldots, n$ with equality when $k=n$.
\end{thm}

There is a known characterization of score sequences of digraphs in general but score is defined in a different way. Avery \cite {avery} defined the \emph{score} of a vertex $x$ in a digraph of order $n$ to be $s(x)=n-1+d^+(x)-d^-(x)$ an characterized score sequences of digraphs.

\begin{thm}[\rm Avery]
	The nondecreasing sequence $\sequence s n$ of nonnegative integers $s_i$ is the score sequences of a digraph of order $n$ if and only if
	$$\sum_{i=1}^ks_i\ge 2\binom k 2$$
	for $k=1,2, \ldots, n$ with equality when $k=n$.
\end{thm}

An oriented complete bipartite graph is called a \emph{bitournament}. Moon \cite{moon} characterized the score sequences of a bipartite tournament. Score, as used here is in the usual sense, outdegree.  We will not make use of the notion of score in the sense of Avery.

\begin{thm}[Moon] Two nondecreasing sequences $\sequence a m$ and $\sequence b n$ of nonnegative integers is a pair of score sequences of a bitournament if and only if 
	$$\sum_{i=1}^k a_i+\sum_{j=1}^\ell b_j\ge k\ell$$
	for $1\le k\le m, 1\le \ell\le n$ with equality when $k=m$ and $\ell=n$.
	\label{thm:moon}
\end{thm}

We are going to give another characterization of the score sequences of bitournaments under the usual notion of score.

\section{Preliminaries}
We will formally define here concepts that are needed in our discussion.

\begin{defn}
	Let $X=\{x_1, \ldots, x_m\}$ and $Y=\{y_1, \ldots, y_n\}$ be the partite sets of a bitournament $T$.  Let $a_i$ be the score of $x_i$ and let $b_i$ be the score of $y_i$, and consider the sequences $A=\sequence a m$, $B=\sequence b n$.  We call the pair $\langle A, B\rangle$ \emph{a score sequence} of the bitournament $T$.
\end{defn}

Note that in the definition, we are not requiring a nondecreasing arrangement of the scores.  If the scores are arranged in nondecreasing order, then we call the pair of sequences \emph{the score sequence} of the bitournament.

\begin{defn}
	A sequence $A=\sequence \alpha m$ is called an $(m,n)$-sequence if $0\le \alpha_i\le n$ for $i=1,2,\ldots ,m$ and where each $\alpha_i$ is a real number.
\end{defn}

Observe that if $A$ is an $(m,n)$-sequence, then $A$ is an $(m,n')$-sequence for $n'>n$.

\begin{defn}
	The \emph{conjugate} of an $(m,n)$-sequence $A=\sequence \alpha n$ is the sequence $\bar A=\langle n-\alpha_i\rangle_{i=1}^m$.
\end{defn}

Note that the conjugate of an $(m,n)$-sequence is also an $(m,n)$-sequence.

\begin{defn}
	A pair $P=\langle A,B\rangle$ of sequences where $A$ is an $(m,n)$-sequence and $B$ is an $(n,m)$-sequence is called an \emph{$(m,n)$-feasible pair} if $\sum_{\alpha\in A}\alpha+\sum_{\beta\in B}\beta=mn$.
\end{defn}

\begin{exmp}
	Let $A=\left\langle 0, \tfrac{1}{3}, 4\right\rangle$, and $B=\left\langle  \tfrac{5}{3},1, 2,3\right\rangle$. Then $P=\langle A, B\rangle$ is a $(3,4)$-feasible pair.
\end{exmp}

Although the elements of sequences are real numbers in general, we shall focus on sequences with integer elements only.\\

A concept that plays a key role in this study is that of $c$-trimming of an $(m,n)$-sequence.

\begin{defn}
	Let $A=\sequence a m$ be an $(m,n)$-sequence, where each $a_i$ is an integer.  If $c>0$ is a constant which is at most equal to the number of positive elements of $A$, then any sequence $A^{\langle c\rangle}$ obtained by subtracting 1 from $c$ positive elements of $A$ is called a \emph{$c$-trimming of $A$}.  If 1 is subtracted from $c$ largest positive elements of $A$, we call  the resulting sequence $A_{\langle c\rangle}$ a \emph{normal $c$-trimming of $A$}.
\end{defn}
 Note that $c$-trimming an $(m,n)$-sequence results into another $(m,n)$-sequence.
 
\begin{exmp} Consider the $(6,6)$-sequence $A=\langle 0, 2, 1, 5, 3, 2\rangle$.  A 3-trimming of $A$ is $A^{\langle 3\rangle}=\langle 0,1,0,4,3,2\rangle$. A normal 3-trimming of $A$ is $A^{\langle 3\rangle}=\langle 0,1,1,4,2,2\rangle$.
\end{exmp}
	
Let us apply a 5-trimming to $A_{\langle 3\rangle}=\langle 0,1,1,4,2,2\rangle$ to obtain $\langle 0,0,0,3,1,1\rangle$.  We shall denote this by $A_{\langle 3,5\rangle}$.  Note that a 5-trimming is not applicable to $A^{\langle 3\rangle}=\langle 0,1,0,4,3,2\rangle$ because it has only 4 positive elements. 

In general, the number of positive elements of a normal $c$-trimming of an $(m,n)$-sequence $A$ is greater than or equal to the number of positive elements of any $c$-trimming of $A$.

For brevity, we shall define 0-trimming by $A_{\langle 0\rangle}=A$. The symbol $A_{\langle c_1, c_2, \ldots, c_k\rangle}$  means a normal $c_1$-trimming, followed by a normal $c_2$-trimming, \emph{etc.}, and lastly a normal  $c_k$-trimming. Note that $A_{\langle c_1, c_2, \ldots, c_k\rangle}$ is not necessarily equal to $A^{\langle c_1, c_2, \ldots, c_k\rangle}$. In fact, if $A_{\langle c_1, c_2,\ldots, c_k\rangle}$ exists, then $A^{\langle c_1, c_2,\ldots, c_k\rangle}$ may not even exist.

For simplicity, if $C=\langle c_1, c_2, \ldots, c_k\rangle$, we shall write $A_C$ for $A_{\langle c_1, c_2, \ldots, c_k\rangle}$ and $A^C$ for $A^{\langle c_1, c_2, \ldots, c_k\rangle}$.

Let us note that if $A^C$ exists, then $A_C$ exists, but not conversely.

\begin{thm}
	Let $A$ be an $(m,n)$-sequence of integers and $B$ an $(n,m)$-sequence of integers such that $\sum_{a\in A}a=\sum_{b\in B}b$.  If $A_B$ exists, then $A_B=\mathbf{0}$. 
\end{thm}

\proof Let $S=\sum_{a\in A}a=\sum_{b\in B}b$. Assume that $A_B$ exists. Let $B=\langle b_1, b_2, \ldots, b_n\rangle$.  Then $\sum_{x\in A_{\langle b_1\rangle}}x=S-b_1$, $\sum_{x\in A_{\langle b_1, b_2\rangle}}=S-b_1-b_2$, and so on.  Therefore, $\sum_{x\in A_B}x=0$. It follows that $A_B=\mathbf{0}$. \qed

\begin{thm}
	Let $P=\langle A, B\rangle$ be an $(m,n)$-feasible pair. Then $\sum_{\bar \beta\in \bar B}\bar \beta=\sum_{\alpha\in A}\alpha$.
\end{thm}

\proof The proof is straightforward.
\begin{align*}
	\sum_{\alpha\in A}\alpha-\sum_{\bar\beta\in \bar B}\bar\beta&=\sum_{\alpha\in A}\alpha-\sum_{\beta\in B}(m-\beta)\\
	&=\sum_{\alpha\in A}\alpha-mn+\sum_{\beta\in B}\beta\\
	&=0
\end{align*}
This proves the Theorem. \qed

\begin{exmp}
	Let $A=\langle 5,3,2, 0\rangle$ and  $B=\langle 1,2,2,2,3\rangle$.  Then $\langle A, B\rangle$ is a $(5,4)$-feasible pair. The conjugate of $B$ is $\bar B=\langle 3,2,2,2,1\rangle$.
	
	We have  $\bar B=\langle 1,2,2,2,3,\rangle$.  Then
	\begin{align*}
		\bar B&=\langle 3,2,2,2,1\rangle\\
		\bar B_{\langle 5\rangle}&=\langle 0,1,1,1,2\rangle\\
		\bar B_{\langle 5,3\rangle}&=\langle 0,1,0,0,1\rangle\\
		\bar B_{\langle 5,3,2\rangle}&=\langle 0,0,0,0,0\rangle\\
		\bar B_A&=\langle 0,0,0,0,0\rangle
	\end{align*}
\end{exmp}

The next result is easy and the proof is omitted.

\begin{thm}
	Let $A$ be an $(m,n)$-sequence, and $B$ an $(n,m)$-sequence such that $B_A=\mathbf{0}$.  If $A'$ is any sequence obtained from $A$ by permuting its elements, then $B_{A'}=\mathbf{0}$
\end{thm}

\begin{exmp}
Let $A=\langle 5,3,2,0\rangle$ and $B=\langle 3,2,2,2,1\rangle$. From a previous example, $B_A=\mathbf{0}$.Let $A'=\langle 2,0,3,5\rangle$, a permutation of $A$.
 Then
\begin{align*}
	B&=\langle 3,2,2,2,1\rangle\\
	B_{\langle 2\rangle}&=\langle 2,1,2,2,1\rangle\\
	B_{\langle 2, 0\rangle}&=\langle 2,1,2,2,1\rangle\\
	B_{\langle 2,0,3\rangle}&=\langle 1,1,1,1,1\rangle\\
	B_{A'}&=\langle 0,0,0,0,0\rangle
\end{align*}

\end{exmp}

For a given sequence $A$ of real numbers, we denote by $cA$ the sequence obtained from $A$ by multiplying every element of $A$ by $c$.
\begin{thm}
	Let $\langle A,B\rangle$ be an $(m,n)$-feasible pair.  Then for every integer $c>0$, $\langle A', B'\rangle$ is a $(cm,cn)$-feasible pair, where $A'=\langle\overbrace{cA,cA,\ldots,cA}^c\rangle$ and $B'=\langle \overbrace{cB,cB,\ldots,cB}^c\rangle$. \label{thm:multiple}
\end{thm}

\proof Let $\langle A,B\rangle$ be an $(m,n)$-feasible pair.  Define $A'=\langle \overbrace{cA,cA,\ldots, cA}^c\rangle$ and $B'=\langle \overbrace{cB,cB,\ldots,cB}^c\rangle$.
\begin{align*}
	\sum_{\alpha\in A'}\alpha&=c^2\sum_{\alpha\in A}\alpha\\
	\sum_{\beta\in B'}\beta&=c^2\sum_{\beta\in B}\beta\\
	\sum_{\alpha\in A'}\alpha+\sum_{\beta\in B'}\beta&=c^2\left[\sum_{\alpha\in A}\alpha+\sum_{\beta\in B}\beta\right]\\
	&=c^2mn\\
	&=(cm)(cn)
\end{align*}
Therefore, $\langle A',B'\rangle$ is a $(cm,cn)$-feasible pair. \qed
~\\

In Theorem \ref{thm:multiple}, the elements of $A$ and $B$ are real numbers in general.  If the elements of $A$ and $B$ are rational numbers, then multiplication by some integer $c>0$ will give us sequences with integer elements.

\section{Main Result}
Moon \ref{thm:moon} has characterized the score sequences of bitournaments. The main result of this paper is a characterization of score sequences of bitournaments utilizing the process of trimming of sequences.

For convenience, let us denote by $\mathbf{0}$ any sequence whose elements are all equal to zero.

\begin{thm}
	Let $A=\sequence a m$ and $B=\sequence b n$ be sequences. The pair $P=\langle A, B\rangle$ is a score sequence of a bitournament if and only if $P$ is an $(m,n)$-feasible pair and $\bar B_A=\mathbf{0}$.
	\label{thm:main}
\label{thm:characterize}
\end{thm}

\proof  First let us assume that $P=\langle A, B\rangle$ is a pair of score sequences of a bitournament $G$. Then $P$ is an $(m,n)$-feasible pair since $\sum_{a\in A}a+\sum_{b\in B}b$ is the total number of arcs of $G$, which is equal to $mn$.
	
Let $x_i$ be the vertex with score $a_i$ and let $y_j$ be the vertex with score $b_j$.  Then $\bar B$ is a sequence of indegrees of the vertices $y_j$.  Specifically, the indegree of $y_j$ is $\bar b_j=n-b_j$.

Let us start drawing the oriented complete bipartite graph (bitournament). Start with $x_1$.  Draw all the arcs of $G$ with initial vertex at $x_1$. Then the $a_1$ outgoing arcs land at some vertices $y_j$.  Subtract 1 from each element of $\bar B$ corresponding to these $a_1$ vertices.  The result is $\bar B^{\langle a_1\rangle}$. If we do this up to $a_m$, we would get $\bar B^A=\mathbf{0}$ because we would have accounted for all the arcs  having $y_j$ as the terminal vertex. Note that symmetrically, we can show tht $\bar A_B=\mathbf{0}$.

Conversely, assume that $A=\sequence a m$ is an $(m,n)$-sequence, $B=\sequence b n$ is an $(n,m)$-sequence such that $\sum_{i=1}^ma_i+\sum_{j=1}^n b_i=mn$, and $\bar B_A=\mathbf{0}$. If $A=\mathbf{0}$, then $B=\langle m,m, \ldots, m\rangle$  and $P$ is obviously a pair of score sequences of a bitournament. Likewise, if $B=\mathbf{0}$, then $A=\langle n,n,\ldots, n\rangle$ and $P$ is a pair of score sequences of a bitournament. So let us assume that $A\ne \textbf{0}$, $B\ne\textbf{0}$.

Consider $\bar B=\langle m-b_1, m-b_2, \ldots, m-b_n\rangle$.  We assume that $\bar B_A=\textbf{0}$. Since $\sum_{i=1}^ma_i+\sum_{j=1}^nb_j=mn$, it follows that $\sum_{i=1}^ma_i=\sum_{j=1}^n\bar b_j$. We are now going to construct an oriented complete bipartite graph with $m$ vertices $x_i$ such that the score of $x_i$ is $a_i$; and $n$ vertices $y_j$ such that the score of $y_j$ is $\bar b_j=n-b_j$. 

From the vertex $x_1$ draw $a_i$ arcs with terminal vertices at the $a_1$ vertices $y_j$ corresponding to $a_i$ largest elements of $\bar B$. Then subtract 1 from the outdegree of each vertex involved. This changes $\bar B$ to $\bar B_{\langle a_1\rangle}$.  We go next to $x_2$ and do a similar thing.  Continue until we reach $x_m$. At this stage, we have made the score of $x_i$ equal to $a_i$, for $1\le i\le m$.

To complete our oriented graph, we draw all missing arcs and orient them from a $y$-vertex to an $x$-vertex.  This will make the score of $y_j$ equal to $m-\bar b_j=b_j$.  Thus, the result is a bipartite tournament with a pair of score sequences $\langle A, B\rangle$. \qed

\begin{exmp}
	Determine whether or not the pair $\langle A, B\rangle$ is the score sequence of a bitournament, where $A=\langle 1,1,2,2,3,4\rangle$, and $B=\langle 1,2,3,5,6\rangle$.
	
	We see that $A$ is a $(6,5)$-sequence and $B$ is a $(5,6)$-sequence, and $\sum_{i=1}^6a_i+\sum_{j=1}^5b_j=13+17=6\times 5$. The conjugate of $B$ is $\bar B=\langle 5,4,3,1,0\rangle$. 
	\begin{align*}
		\bar B&=\langle5,4,3,1\rangle\\
		\bar B_{\langle1\rangle}&=\langle 4,4,3,1\rangle\\
		\bar B_{\langle 1,1\rangle}&=\langle 3,4,3,1\rangle\\
		\bar B_{\langle 1,1,2\rangle}&=\langle2,3,3,1\rangle\\
		\bar B_{\langle 1,1,2,2\rangle}&=\langle 2,2,2,1\rangle\\
		\bar B_{1,1,2,2,3\rangle}&=\langle 1,1,1,1\rangle\\
		\bar B_{\langle 1,1,2,2,3,4\rangle}&=\langle 0,0,0,0\rangle
	\end{align*}
Therefore, $\langle A, B\rangle$ is a pair of score sequences for a bipartite tournent.  Following the constructive proof of Theorem \ref{thm:main}, let us draw the bipartite tournament. The step-by-step construction is illustrated below.
\end{exmp}
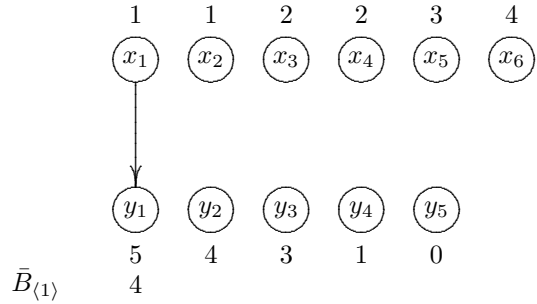
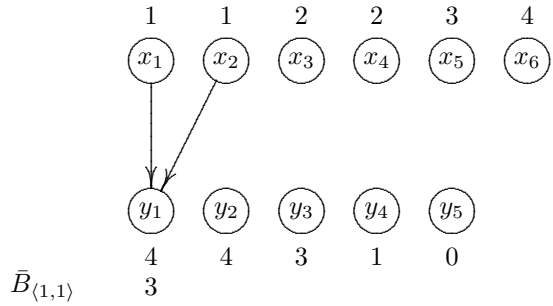
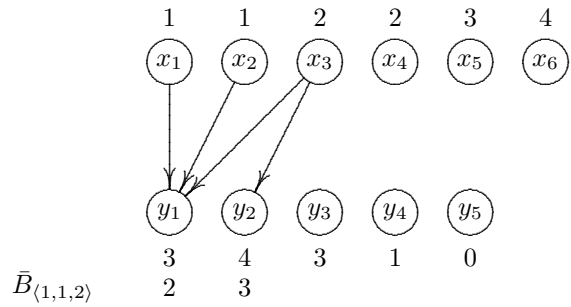
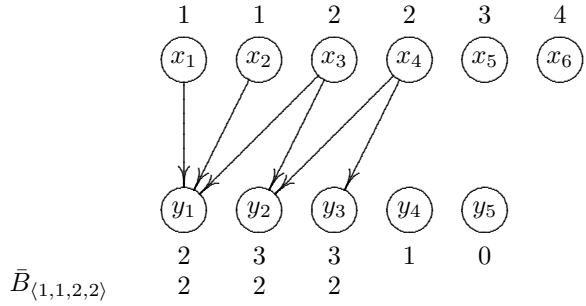
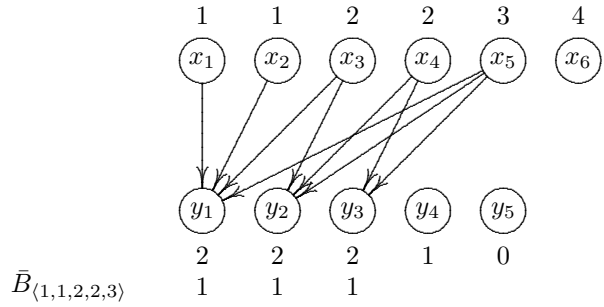
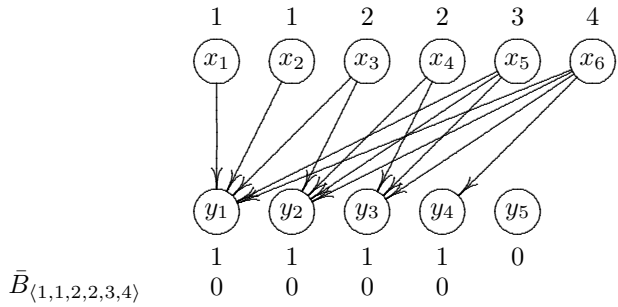
\newbox\btournament
\setbox\btournament=\hbox{$\pic
	\xunit1cm
	\yunit1cm
	\VertexRadius0.3cm
	\Vertex(1,0) (2,0) (3,0) (4,0) (5,0) 
	\Vertex(1,2) (2,2) (3,2) (4,2) (5,2) (6,2)
	\Align[c] ($x_1$) (1,2)
	\Align[c] ($x_2$) (2,2)
	\Align[c] ($x_3$) (3,2)
	\Align[c] ($x_4$) (4,2)
	\Align[c] ($x_5$) (5,2)
	\Align[c] ($x_6$) (6,2)
	\Align[c] ($y_1$) (1,0)
	\Align[c] ($y_2$) (2,0)
	\Align[c] ($y_3$) (3,0)
	\Align[c] ($y_4$) (4,0)
	\Align[c] ($y_5$) (5,0)
	\cip$}
	
$$\pic
\xunit1cm
\yunit1cm
\VertexRadius0.3cm
\Align[ll] (\copy\btournament) (0.7,-0.3)
\Align[c] (1) (1,2.6) (2,2.6 ) (4,-0.6)
\Align[c] (2) (3,2.6) (4,2.6) 
\Align[c] (3) (5,2.6) (3,-0.6)
\Align[c] (4) (6,2.6) (2,-0.6)
\Align[c] (5) (1,-0.6)
\Align[c] (0) (5,-0.6)
\Align[r] ($\Bar B_{\langle 1\rangle}$) (0,-1)
\Align[c] (4) (1,-1)
\Arc(1,2) (1,0)
\cip$$
	
$$\pic
\xunit1cm
\yunit1cm
\VertexRadius0.3cm
\Align[ll] (\copy\btournament) (0.7,-0.3)
\Align[c] (1) (1,2.6) (2,2.6 ) (4,-0.6)
\Align[c] (2) (3,2.6) (4,2.6) 
\Align[c] (3) (5,2.6) (3,-0.6)
\Align[c] (4) (6,2.6) (2,-0.6) (1,-0.6)
\Align[c] (0) (5,-0.6)
\Align[r] ($\Bar B_{\langle 1,1\rangle}$) (0,-1)
\Align[c] (3) (1,-1)
\Arc(1,2) (1,0)
\Arc(2,2) (1,0)
\cip$$
	
$$\pic
\xunit1cm
\yunit1cm
\VertexRadius0.3cm
\Align[ll] (\copy\btournament) (0.7,-0.3)
\Align[c] (1) (1,2.6) (2,2.6 ) (4,-0.6)
\Align[c] (2) (3,2.6) (4,2.6) 
\Align[c] (3) (5,2.6) (3,-0.6)  (1,-0.6) (2,-1)
\Align[c] (4) (6,2.6) (2,-0.6) 
\Align[c] (0) (5,-0.6)
\Align[r] ($\Bar B_{\langle 1,1,2\rangle}$) (0,-1)
\Align[c] (2) (1,-1)
\Arc(1,2) (1,0)
\Arc(2,2) (1,0)
\Arc(3,2) (1,0) (3,2) (2,0)
\cip$$	

$$\pic
\xunit1cm
\yunit1cm
\VertexRadius0.3cm
\Align[ll] (\copy\btournament) (0.7,-0.3)
\Align[c] (1) (1,2.6) (2,2.6 ) (4,-0.6)
\Align[c] (2) (3,2.6) (4,2.6) (1,-0.6)
\Align[c] (3) (5,2.6) (3,-0.6)  (2,-0.6)
\Align[c] (4) (6,2.6)  
\Align[c] (0) (5,-0.6)
\Align[r] ($\Bar B_{\langle 1,1,2,2\rangle}$) (0,-1)
\Align[c] (2) (1,-1) (2,-1) (3,-1)
\Arc(1,2) (1,0)
\Arc(2,2) (1,0)
\Arc(3,2) (1,0) (3,2) (2,0)
\Arc(4,2) (2,0) (4,2) (3,0)
\cip$$	

$$\pic
\xunit1cm
\yunit1cm
\VertexRadius0.3cm
\Align[ll] (\copy\btournament) (0.7,-0.3)
\Align[c] (1) (1,2.6) (2,2.6 ) (4,-0.6)
\Align[c] (2) (3,2.6) (4,2.6)  (3,-0.6)  (1,-0.6)  (2,-0.6)
\Align[c] (3) (5,2.6)
\Align[c] (4) (6,2.6)  
\Align[c] (0) (5,-0.6)
\Align[r] ($\Bar B_{\langle 1,1,2,2,3\rangle}$) (0,-1)
\Align[c] (1) (1,-1) (2,-1) (3,-1)
\Arc(1,2) (1,0)
\Arc(2,2) (1,0)
\Arc(3,2) (1,0) (3,2) (2,0)
\Arc(4,2) (2,0) (4,2) (3,0)
\Arc(5,2) (1,0) (5,2) (2,0) (5,2) (3,0)
\cip$$	

$$\pic
\xunit1cm
\yunit1cm
\VertexRadius0.3cm
\Align[ll] (\copy\btournament) (0.7,-0.3)
\Align[c] (1) (1,2.6) (2,2.6 ) (3,-0.6)  (1,-0.6)  (2,-0.6) (4,-0.6)
\Align[c] (2) (3,2.6) (4,2.6) 
\Align[c] (3) (5,2.6)
\Align[c] (4) (6,2.6)  
\Align[c] (0) (5,-0.6) (4,-1)
\Align[r] ($\Bar B_{\langle 1,1,2,2,3,4\rangle}$) (0,-1)
\Align[c] (0) (1,-1) (2,-1) (3,-1)
\Arc(1,2) (1,0)
\Arc(2,2) (1,0)
\Arc(3,2) (1,0) (3,2) (2,0)
\Arc(4,2) (2,0) (4,2) (3,0)
\Arc(5,2) (1,0) (5,2) (2,0) (5,2) (3,0)
\Arc(6,2) (1,0) (6,2) (2,0) (6,2) (3,0) (6,2) (4,0)
\cip$$	

After these, we draw upward arcs between pairs of vertices that are not yet connected by an arc.

\begin{thmcor}
	Let $P=\langle A, B\rangle$ be an $(m,n)$-feasible pair.  Then $\bar A_B=\mathbf{0}$ if and only if $\bar B_A=\mathbf{0}$.
\end{thmcor}

\section{Concluding Remarks}
Two examples will be given here for the purpose of illustrating and comparing Moon's theorem and Theorem \ref{thm:main}, the main result of this paper.

\begin{exmp}
	Let $A=\langle 1,3,4,5\rangle$ and $B=\langle 0,1,2,2,2\rangle$.  Show that $\langle A, B\rangle$ is a pair of score sequences for some bitournament.\\
	
 We will use first Theorem \ref{thm:main} to do the verification. Here, $A$ is a $(4,5)$-sequence, $B$ is a $(5,4)$-sequence and the sum of all their elements is $20=4\times 5$.		
		\begin{align*}
			\bar B&=\langle 4,3,2,2,2\rangle\\
			\bar B_{\langle 1\rangle}&=\langle3,3,2,2,2\rangle\\
			\bar B_{\langle 1,3\rangle}&=\langle 2,2,1,2,2\rangle\\
			\bar B_{\langle 1,3,4\rangle}&=\langle1,1,1,1,1\rangle\\
			\bar B_{\langle 1,3,4,5\rangle}&=\langle 0,0,0,0,0\rangle
		\end{align*}
		Therefore, $\langle A,B\rangle$ is the score sequence of some bitournament with underlying graph $K_{4,5}$, a complete bipartite graph. 
		
		Let us use Moon's Theorem \ref{thm:moon} this time. We tabulate, for each pair of values $(k,\ell)$, the sum $\sum_{i=1}^ka_i+\sum_{j=1}^\ell b_i$ and the product $k\ell$.  We see that the sum is not less than the product $k\ell$ in every case, and that equality holds when $k=4$ and $\ell=5$.
			
			$$\begin{tabular}{|c|c|c|c|c|c|c|c|c|c|c|c|c|c|c|c|c|c|c|c|c|}
				\hline
				$k$&1&1&1&1&1&2&2&2&2&2&3&3&3&3&3&4&4&4&4&4\\
				\hline
				$\ell$&1&2&3&4&5&1&2&3&4&5&1&2&3&4&5&1&2&3&4&5\\
				\hline
				sum&1&2&4&6&8&4&5&7&9&11&8&9&11&13&15&13&14&16&18&20\\
				\hline
				$k\ell$&1&2&3&4&5&2&4&6&8&10&3&6&9&12&15&4&8&12&16&20\\
				\hline
			\end{tabular}$$
Therefore, $\langle A, B\rangle$ ia a sore sequence of some bitournament.		
\end{exmp}

\begin{exmp}
	Let $p$ be any positive integer and test the pair $\langle A, B\rangle$, where $A=\langle\overbrace{p,\ldots, p}^{2p}\rangle=B$.
	
	By Moon's theorem, $\sum_{i=1}^kp+\sum_{j=1}^\ell p=kp+\ell p\ge k\tfrac{\ell}{2}+\ell\tfrac{k}{2}=kl$, and when $k=2p=\ell$, the sum is $2p^2+2p^2=4p^2=(2p)(2p)=mn$. 
	
	By Theorem \ref{thm:main}, $\langle A, B\rangle$ is a $(2p,2p)$-feasible pair. The conjugate of $B$ is $\bar B=B$, and so  $\bar B_A=\mathbf{0}$.
	
	Both tests confirm that $\langle A, B\rangle$ is a pair of score sequences of some bitournament.
\end{exmp}

\end{document}